\documentclass[12pt]{article}
\usepackage{amssymb,amsfonts,amsmath, psfrag,eepic,colordvi,graphicx,epsfig,ytableau}
\usepackage{amssymb,latexsym,graphics}
\usepackage{MnSymbol}
\usepackage[enableskew]{youngtab}
\usepackage{float}
\usepackage{tikz}
\topskip  
\newcommand{\rmnum}[1]{\romannumeral #1}

 \numberwithin{equation}{section}
\newtheorem{theorem}{Theorem}[section]

\newtheorem{corollary}[theorem]{Corollary}

\newtheorem{lemma}[theorem]{Lemma}

\newtheorem{observation}[theorem]{Observation}

\begin{document}
\parskip 6pt

\pagenumbering{arabic}
\def\sof{\hfill\rule{2mm}{2mm}}
\def\ls{\leq}
\def\gs{\geq}
\def\SS{\mathcal S}
\def\qq{{\bold q}}
\def\MM{\mathcal M}
\def\TT{\mathcal T}
\def\EE{\mathcal E}
\def\lsp{\mbox{lsp}}
\def\rsp{\mbox{rsp}}
\def\pf{\noindent {\it Proof.} }
\def\mp{\mbox{pyramid}}
\def\mb{\mbox{block}}
\def\mc{\mbox{cross}}
\def\qed{\hfill \rule{4pt}{7pt}}
\def\block{\hfill \rule{5pt}{5pt}}

\begin{center}
{\Large\bf    Quasi-Stirling Permutations on Multisets  }
\end{center}

\begin{center}
{\small Sherry H.F. Yan,  Lihong Yang, Yunwei Huang, Xue Zhu }

 Department of Mathematics\\
Zhejiang Normal University\\
 Jinhua 321004, P.R. China

 hfy@zjnu.cn

\end{center}

\noindent {\bf Abstract.} A permutation $\pi$ of a multiset is said to be a {\em quasi-Stirling
} permutation  if  there does not exist four indices $i<j<k<\ell$ such that $\pi_i=\pi_k$ and $\pi_j=\pi_{\ell}$. Define $$ \overline{Q}_{\mathcal{M}}(t,u,v)=\sum_{\pi\in \overline{\mathcal{Q}}_{\mathcal{M}}}t^{des(\pi)}u^{asc(\pi)}v^{plat(\pi)},$$ where $\overline{\mathcal{Q}}_{\mathcal{M}}$ denotes the set of quasi-Stirling
 permutations on the multiset $\mathcal{M}$,  and  $asc(\pi)$ (resp. $des(\pi)$, $plat(\pi)$) denotes  the number of   ascents (resp. descents, plateaux) of $\pi$.  Denote by $\mathcal{M}^{\sigma}$ the multiset   $\{1^{\sigma_1}, 2^{\sigma_2}, \ldots, n^{\sigma_n}\}$, where $\sigma=(\sigma_1, \sigma_2, \ldots, \sigma_n)$ is an $n$-composition of $K$ for positive integers $K$ and $n$.  In this paper, we show that  $\overline{Q}_{\mathcal{M}^{\sigma}}(t,u,v)=\overline{Q}_{\mathcal{M}^{\tau}}(t,u,v)$ for any two $n$-compositions $\sigma$ and $\tau$  of $K$. This is accomplished by establishing an $(asc, des, plat)$-preserving bijection
    between  $\overline{\mathcal{Q}}_{\mathcal{M}^{\sigma}}$ and $\overline{\mathcal{Q}}_{\mathcal{M}^{\tau}}$.  As applications,  we obtain  generalizations of  several results  for quasi-Stirling permutations on $\mathcal{M}=\{1^k,2^k, \ldots, n^k\}$ obtained by Elizalde and solve an open problem posed by Elizalde.




\section{Introduction}

Let $\mathcal{M}=\{1^{k_1}, 2^{k_2}, \ldots, n^{k_n}\}$ be a multiset where $k_i$ is the number of occurrences of $i$ in  $\mathcal{M}$ and $k_i\geq 1$.
  A permutation $\pi$ of a multiset $\mathcal{M}$ is said to be  a {\em Stirling } permutation if   $i<j<k$ and $\pi_i=\pi_k$, then $\pi_j>\pi_i$.
  Stirling permutations were originally introduced by Gessel and Stanley \cite{Gessel} in the case of the multiset $\mathcal{M}=\{1^{2}, 2^{2}, \ldots, n^{2}\}$. Park \cite{Park, Park2}  and Liu \cite{Liu} studied  the distribution of various statistics over Stirling permutations on $\mathcal{M}=\{1^{k}, 2^{k}, \ldots, n^{k}\}$. Brenti \cite{Brenti1, Brenti2} and Dzhumadil'daev and  Yeliussizov \cite{Dzhumadil'daev} investigated the   descent polynomials of Stirling permutations on an arbitrary multiset.
In   analogy to Stirling permutations, Archer, Gregory, Pennington and Slayden \cite{Archer} introduced   {\em quasi-Stirling} permutations. A permutation $\pi$ of a multiset is said to be a {\em quasi-Stirling
} permutation  if  there does not exist four indices $i<j<k<\ell$ such that $\pi_i=\pi_k$ and $\pi_j=\pi_\ell$. For a multiset $\mathcal{M}$, denote by   $\overline{\mathcal{Q}}_{\mathcal{M}}$ the set of quasi-Stirling permutations of $\mathcal{M}$. For example, if $\mathcal{M}=\{1^2, 2^2\}$, we have
$
\overline{\mathcal{Q}}_{\mathcal{M}} =\{1221, 2112, 1122, 2211\}.
$

For a permutation $\pi=\pi_1\pi_2\ldots \pi_n$,   an index $i$,  $0\leq i\leq n$, is called an  {\em ascent} (resp. {\em a descent},  {\em a plateau}) of $\pi$ if $\pi_{i}<\pi_{i+1}$ (resp. $\pi_i>\pi_{i+1}$, $\pi_i=\pi_{i+1}$) with the convention $\pi_0=\pi_{n+1}=0$.  Let $asc(\pi)$ (resp. $des(\pi)$, $plat(\pi)$) denote the number of   ascents (resp. descents, plateaux) of $\pi$. In \cite{Bona}, Bona proved that   the three
statistics $asc$, $plat$ and $des$ are equidistributed over Stirling permutations on $\mathcal{M}=\{1^2, 2^2, \ldots, n^2\}$.
Define $$ \overline{Q}_{\mathcal{M}}(t,u,v)=\sum_{\pi\in \overline{\mathcal{Q}}_{\mathcal{M}}}t^{des(\pi)}u^{asc(\pi)}v^{plat(\pi)}.$$
 The polynomial  $ \overline{Q}_{\mathcal{M}}(t,1,1)$  is called
the  {\em quasi-Stirling polynomial } on the multiset $\mathcal{M}$.
 By employing generating function arguments,  Elizalde \cite{Elizalde}     derived that
  \begin{equation}\label{eq3}
 \sum_{m=0}^{\infty} {m^n\over n+1}{n+m\choose m}t^m={\overline{Q}_{\mathcal{M}}(t,1,1)\over (1-t)^{2n+1}},
 \end{equation}
when $\mathcal{M}=\{1^2, 2^2, \ldots, n^2\}$
   and asked for a combinatorial  proof.
Recently, Yan and Zhu \cite{Yanzhu} provided such a combinatorial proof and  further derived that
\begin{equation}\label{eqmain}
\sum_{m=0}^{\infty} {m^n\over K-n+1}{K-n+m \choose m} t^m={\overline{Q}_{\mathcal{M}}(t,1,1)\over (1-t)^{K+1}},
\end{equation}
where $\mathcal{M}=\{1^{k_1}, 2^{k_2}, \ldots, n^{k_n}\}$ and $K=k_1+k_2+\ldots+k_n$ with $k_i\geq 1$.

For positive integers $K$ and $n$, let $\sigma=(\sigma_1, \sigma_2, \ldots, \sigma_n)$ be an $n$-composition of $K$. Denote by $\mathcal{M}^{\sigma}$ the multiset   $\{1^{\sigma_1}, 2^{\sigma_2}, \ldots, n^{\sigma_n}\}$.
By (\ref{eqmain}), for any two $n$-compositions $\sigma$ and $\tau$  of $K$,  we have   $\overline{Q}_{\mathcal{M}^{\sigma}}(t,1,1)=\overline{Q}_{\mathcal{M}^{\tau}}(t,1,1)$.   In this paper, we get the following generalizations.

\begin{theorem}\label{mainth1}
Let $\mathcal{M}=\{1^{k_1}, 2^{k_2}, \ldots, n^{k_n}\}$, $\mathcal{M'}=\{1^{K-n+1}, 2,3,\ldots, n\}$ and $K=k_1+k_2+\ldots+k_n$ with $k_i\geq 1$.
There is a bijection $\Phi$ between $\overline{\mathcal{Q}}_{\mathcal{M}}$ and $\overline{\mathcal{Q}}_{\mathcal{M'}}$ such that for any $\pi\in \overline{\mathcal{Q}}_{\mathcal{M}}$, we have $$(asc, des, plat)\pi= (asc, des, plat)\Phi(\pi).$$
Consequently, we have
      $\overline{Q}_{\mathcal{M}}(t,u,v)=\overline{Q}_{\mathcal{M'}}(t,u,v)$.
    \end{theorem}

The following result follows immediately from Theorem \ref{mainth1}.
\begin{theorem}\label{mainth}
Let  $K$ and $n$ be positive integers. For any two $n$-compositions $\sigma$ and $\tau$  of $K$,  there exists a bijection  between  $\overline{\mathcal{Q}}_{\mathcal{M}^{\sigma}}$ and $\overline{\mathcal{Q}}_{\mathcal{M}^{\tau}}$  which preserves   the number of ascents, descents and plateaux.
Consequently, we have
      $\overline{Q}_{\mathcal{M}^{\sigma}}(t,u,v)=\overline{Q}_{\mathcal{M}^{\tau}}(t,u,v)$.
\end{theorem}

     Let $\mathcal{J}_{n,r}$ be the set of permutations $\sigma=\sigma_1\sigma_2\ldots\sigma_{n-r}$ where each $\sigma_i\in [n]$ and $\sigma_i\neq \sigma_j$ for all $1\leq i,j\leq n-r$. Recall that for a permutation $\pi=\pi_1\pi_2\ldots\pi_n$, an index $i$ is said to be an {\em excedance} of $\pi$ if $\pi_i>i$. Let $exc(\pi)$ denote the number of excedances of $\pi$. By employing generating function arguments, Elizalde \cite{Elizalde} derived that for $\mathcal{M}=\{1^2, 2^2, \ldots, n^2\}$,   quasi-Stirling permutations $\pi\in \overline{\mathcal{Q}}_{\mathcal{M}}$ with $des(\pi)=d+1$  are equinumerous with  permutations $\sigma\in  \mathcal{J}_{2n,n+1}$ with $ exc(\sigma)=d$, and asked for a combinatorial proof.  Relying on Theorem \ref{mainth1}, we present such a combinatorial proof and further extend Elizalde's result to an arbitrary multiset.

    \begin{theorem}\label{mainth2}
    Let $\mathcal{M}=\{1^{k_1}, 2^{k_2}, \ldots, n^{k_n}\}$ and $K=k_1+k_2+\ldots+k_n$ with $k_i\geq 1$.
There is a bijection between  quasi-Stirling permutations $\pi\in \overline{\mathcal{Q}}_{\mathcal{M}}$ with $des(\pi)=d+1$ and permutations $\sigma\in  \mathcal{J}_{K,K-n+1}$ with $ exc(\sigma)=d$.
    \end{theorem}

 By establishing a bijection between   quasi-Stirling  permutations  and edge-labeled plane (ordered) rooted trees, Elizalde \cite{Elizalde} confirmed  a   conjecture  posed by Archer et al. \cite{Archer} which asserts that the number of $\pi\in \overline{\mathcal{Q}}_{\mathcal{M}}$ with $des(\pi)=n$ is equal to $(n+1)^{n-1}$ for $\mathcal{M}=\{1^2, 2^2, \ldots, n^2\}$. He further proved that
 the number of $\pi\in \overline{\mathcal{Q}}_{\mathcal{M}}$ with $des(\pi)=n$ is equal to $((k-1)n+1)^{n-1}$ for  $\mathcal{M}=\{1^k, 2^k, \ldots, n^k\}$ and $k\geq 2$.
As an application of Theorem \ref{mainth1}, we obtain the following generalization of Elizalde's result.
\begin{corollary}\label{coro1}
Let $\mathcal{M}=\{1^{k_1}, 2^{k_2}, \ldots, n^{k_n}\}$ and $K=k_1+k_2+\ldots+k_n$ with $k_i\geq 1$. Then the number of permutations  $\pi\in \overline{\mathcal{Q}}_{\mathcal{M}}$ with $des(\pi)=n$ is equal to $(K-n+1)^{n-1}$.
    \end{corollary}

Denote by $\mathcal{S}_n$  the set of all permutations of $[n]$. The polynomial
 $$A_n(t,u)=\sum_{\pi\in \mathcal{S}_n} t^{des(\pi)}u^{asc(\pi)}$$
 is called the {\em   bivariate Eulerian polynomial}.
 Define
 $$
 \widetilde{A}(t,u;z)=1+\sum_{n\geq 1} A_n(t,u){z^n\over n!}.
 $$
In  \cite{Elizalde},   Elizalde  derived that
\begin{equation}\label{eqk}
       \overline{Q}_{\mathcal{M}}(t,u,v)={n!\over (k-1)n+1}[z^n](\widetilde{A}(t,u;z)-1+v)^{(k-1)n+1}.
    \end{equation}
 for $\mathcal{M}=\{1^k,2^k, \ldots, n^k\}$.
 Based on Theorem \ref{mainth1}, we  extend  (\ref{eqk}) to an arbitrary multiset.
 \begin{corollary}\label{coro2}
Let $\mathcal{M}=\{1^{k_1}, 2^{k_2}, \ldots, n^{k_n}\}$ and $K=k_1+k_2+\ldots+k_n$ with $k_i\geq 1$. Then
     $$
     \overline{Q}_{\mathcal{M}}(t,u,v)={n!\over K-n+1}[z^n](\widetilde{A}(t,u;z)-1+v)^{K-n+1}.
     $$
    \end{corollary}

    The rest of this paper is organized as follows. In Section 2,  we shall establish the bijection $\Phi$, thereby  proving  Theorem \ref{mainth1}. As an application of Theorem  \ref{mainth1}, we get the enumeration of quasi-Stirling permutations  with $n$ descents. Section 3 is devoted to the proof of Theorem \ref{mainth2}. In Section 4, relying on Theorem \ref{mainth1}, we express the polynomial  $\overline{Q}_{\mathcal{M}}(t,u,v)$ in terms of the generating function $\widetilde{A}(t,u;z)$, thereby proving Corollary \ref{coro2}.
\section{  Proof of Theorem \ref{mainth1}}

\begin{figure}[h]
\centerline{\includegraphics[width=11cm]{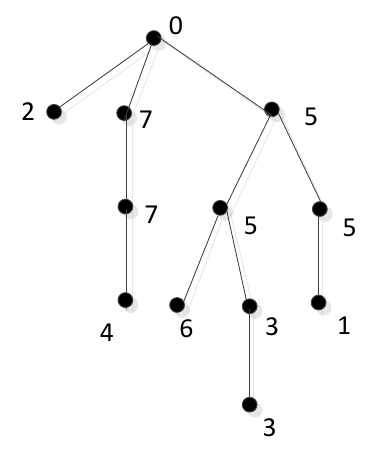}}
 \caption{ A tree $T\in \mathcal{T}_{
\mathcal{M}}$ with $\mathcal{M}=\{1,2,3^{2},4,5^{3},6,7^{2}\}$. }\label{T}
\end{figure}

The objective of this section is  to establish the bijection  $\Phi$ stated  in  Theorem \ref{mainth1}.
Recall that an  {\em ordered} tree is a tree with one designated vertex, which is called the root, and the subtrees of
each vertex are linearly ordered.  In an ordered tree $T$,   the {\em level} of a vertex $v$ in $T$  is defined to be the length of the unique path from the root to $v$.  A vertex $v$ is said to be {\em at odd level} (resp. {\em at even level}) if the level of $v$ is odd (resp. even).   In a rooted tree, the vertex $v$ is a {\em descendant} of the vertex $u$ if $u$ lies on the path from $v$ to the root.
 In order to prove (\ref{eqmain}),  Yan and Zhu \cite{Yanzhu}  introduced   a new
class of ordered   labeled trees $T$   verifying the following properties:
\begin{itemize}
\item[{\upshape (\rmnum{1})}] the vertices are labeled by the elements of the multiset  $\{0\}\cup \mathcal{M}$, where $\mathcal{M}=\{1^{k_1}, 2^{k_2}, \ldots, n^{k_n}\}$ with $k_i\geq 1$;
    \item[{\upshape (\rmnum{2})}] the root is labeled by $0$;
    \item[{\upshape (\rmnum{3})}] for  a vertex  $v$ at odd level,  if $v$ is labeled by $i$, then $v$ has  exactly  $k_i-1$ children  and the children of $v$ have the same label as that of $v$
  in $T$.
\end{itemize}
Let $\mathcal{T}_{\mathcal{M}}$ denote the set of  such ordered    labeled trees. For example,  a tree $T\in \mathcal{T}_{
\mathcal{M}}$ with $\mathcal{M}=\{1,2,3^{2},4,5^{3},6,7^{2}\}$ is illustrated in Figure \ref{T}.

The number of cyclic descents  of a sequence $\pi=\pi_1\pi_2\ldots\pi_n$ is defined  to
be
$$
cdes(\pi)=|\{i\mid \pi_i>\pi_{i+1}\}|
$$
with the convention $\pi_{n+1}=\pi_1$.
Similarly, the number of cyclic ascents  of a sequence $\pi=\pi_1\pi_2\ldots\pi_n$ is defined to be

$$
casc(\pi)=|\{i\mid \pi_i<\pi_{i+1}\}|
$$
with the convention $\pi_{n+1}=\pi_1$.

In an ordered  labeled tree $T$, let $u$ be a vertex of $T$. Suppose that $u$ has $\ell$ children $v_1, v_2, \ldots, v_{\ell}$ listed from left to right. The number of  cyclic descents (resp. cyclic ascents) of $u$, denoted by $cdes(u)$ (resp. $casc(u)$),  is defined to be  the number  $cdes(uv_1v_2\ldots v_{\ell})$ (resp. $casc(uv_1v_2\ldots v_{\ell})$). The number of cyclic descents and cyclic ascents of $T$  are defined to be
$$
cdes(T)=\sum_{u\in V(T)}cdes(u),
$$
and
$$
casc(T)=\sum_{u\in V(T)}casc(u),
$$
where $V(T)$ denotes the vertex set of $T$.
Denote by $eleaf(T)$ the number of leaves at even level.
 For example, if we let $T$ be a tree as shown  in Figure \ref{T}, we have $casc(T)=6$, $cdes(T)=5$ and  $eleaf(T)=1$.
Let $\pi$ be a permutation, the leftmost (resp. rightmost) entry of $\pi$ is denoted by $first(\pi)$ (resp. $last(\pi)$). Similarly, for an ordered labeled tree $T$, we denote by $first(T)$   (resp. $last(T)$) the leftmost (resp. rightmost) child of the root.

In order to prove (\ref{eqmain}),  Yan and Zhu \cite{Yanzhu} established a bijection $\phi$  between $\mathcal{T}_{\mathcal{M}}$ and $
\overline{\mathcal{Q}}_{\mathcal{M}}$. First we give an overview of the bijection $\phi$.
If $T$ has only one vertex, let $\phi(T)=\epsilon$, where $\epsilon$ denotes the empty permutation.  Otherwise, suppose that $first(T)=r$.

\noindent {\bf Case 1.}  The leftmost child of the root is a leaf. Let $T_0$ be the tree obtained from $T$ by removing  the leftmost child of the root together with  the edge incident to it. Define
    $\phi(T)=r\phi(T_0)$.

    \noindent {\bf Case 2.} The leftmost child of the root has $k$ children. For $1\leq i\leq k$,  let $T_i$ be the subtree rooted at the $i$-th child of the leftmost child of the root.
         Denote by $T'_i$ the tree obtained from $T_i$ by relabeling its root by $0$.
         Let $T_0$ be the tree obtained from $T$ by removing all the subtrees $T_1$, $T_2$, $\ldots,$ $T_k$ and  all the vertices labeled by $r$ together with  the edges incident to them. Define $\phi(T)=r\phi(T'_1)r\phi(T'_2)\ldots r\phi(T'_k)r\ \phi(T_0)$.

For instance, let $T$ be a tree illustrated in Figure \ref{T}. By applying the map  $\phi$, we have $\phi(T)=27475633515$.

In \cite{Yanzhu}, Yan and Zhu proved that the bijection $\phi$ has the following property.

\begin{theorem}{\upshape   ( \cite{Yanzhu}) }\label{thphi0}
For any nonempty multiset $\mathcal{M}$,the map $\phi$ is  a bijection between $\mathcal{T}_{\mathcal{M}}$ and $
\overline{\mathcal{Q}}_{\mathcal{M}}$  such that $$(cdes,  first, last)T= (des,  first, last)\phi(T)$$ for any $T\in \mathcal{T}_{\mathcal{M}}$.
\end{theorem}

\begin{theorem}\label{thphi}

For any nonempty multiset $\mathcal{M}$, the map $\phi$ induces a bijection between $\mathcal{T}_{\mathcal{M}}$ and $
\overline{\mathcal{Q}}_{\mathcal{M}}$  such that $$(cdes, casc, eleaf,  first, last)T= (des, asc, plat, first, last)\phi(T)$$ for any $T\in \mathcal{T}_{\mathcal{M}}$.
\end{theorem}

\pf  By theorem \ref{thphi0}, it remains to show that $( casc, eleaf)T= (asc, plat)\phi(T)$. We  shall  prove the statement by induction on the cardinality of the multiset $\mathcal{M}$. Suppose that $first(T)=r$. If $|\mathcal{M}|=1$, then $T$ consists of the root $0$ and its child $r$. According to the construction of $\phi$, we have $\phi(T)=r$. Then it is easy to check that       $casc(T)=asc(\phi(T))=1$ and $eleaf(T)=plat(\phi(T))=0$.  Assume that $(casc, eleaf)T= (asc, plat)\phi(T)$ for any $T\in \mathcal{T}_{\mathcal{M}'}$    with $|\mathcal{M}'|<|\mathcal{M}|$. We proceed to show that $(casc, eleaf)T= (asc, plat)\phi(T)$ for any $T\in \mathcal{T}_{\mathcal{M}}$.
 We have two cases.

\noindent {\bf Case 1.}    The leftmost child of the root is a leaf. In this case,   we have $\phi(T)=r\phi(T_0)$. It is routine to check that
    $$
    casc(T)=casc(T_0)+\chi(r<first(T_0)).
    $$
    Here $\chi(S)=1$ if  the statement $S$ is true, and $\chi(S)=0$ otherwise.
    Then by induction hypothesis, we have
    $$
    \begin{array}{lll}
    casc(T)&=&casc(T_0)+\chi(r<first(T_0))\\
    &=&asc(\phi(T_0))+\chi(r<first(\phi(T_0)))\\
    &=&asc(\phi(T)).
    \end{array}
    $$
Similarly, we have $
    eleaf(T)=eleaf(T_0).
    $ and $plat(\phi(T_0))=plat(\phi(T))$.
Again by induction hypothesis, we have $eleaf(T)=eleaf(T_0)=plat(\phi(T_0))=plat(\phi(T))$.

 \noindent {\bf Case 2.}  The leftmost child of the root has $k$ children. Then  we have $\phi(T)=r\phi(T'_1)r\phi(T'_2)\ldots r\phi(T'_k)r\ \phi(T_0)$.
 Define
 $$
 \mathcal{I}=\{i\mid  |T_i|>1, 1\leq i\leq k\},
 $$
 where $|T_i|$ denotes the number of vertices of $T_i$.

 \noindent {\bf Subcase 2.1.}  $|T_0|>1$.
 It is easily seen that
    $$
    casc(T)= casc(T_0)+\chi(r<first(T_0))+\sum_{i\in \mathcal{I}}\big(casc(T'_i)-1+\chi(r<first(T'_i)\big)+\chi(r>last(T'_i)),    $$
     and $asc(\phi(T))$ is given by
    $$
    asc(\phi(T_0))+\chi(r<first(\phi(T_0)))+\sum_{i\in \mathcal{I}}\big(asc(\phi(T'_i))-1+\chi(r<first(\phi(T'_i))\big)+\chi(r>last(\phi(T'_i))).
    $$
   By induction hypothesis,  it is easy to verify that $casc(T)=asc(\phi(T))$.
   Similarly, we have
   $$
   eleaf(T)=eleaf(T_0)+ |  [k]\setminus \mathcal{I} |+ \sum_{i\in \mathcal{I}}eleaf(T'_i)
   $$
   and
   $$
   plat(\phi(T))= plat(\phi(T_0))+ |  [k]\setminus \mathcal{I} |+ \sum_{i\in \mathcal{I}}plat(\phi(T'_i)).
   $$
   Again by induction hypothesis, we have $eleaf(T)=plat(\phi(T))$.

 \noindent {\bf Subcase 2.2.}  $|T_0|=1$.
 It is easily seen that
    $$
    casc(T)= 1+\sum_{i\in \mathcal{I}}\big(casc(T'_i)-1+\chi(r<first(T'_i)\big)+\chi(r>last(T'_i)),    $$
     and $asc(\phi(T))$ is given by
    $$
    1+\sum_{i\in \mathcal{I}}\big(asc(\phi(T'_i))-1+\chi(r<first(\phi(T'_i))\big)+\chi(r>last(\phi(T'_i))).
    $$
   By induction hypothesis,  it is easy to verify that $casc(T)=asc(\phi(T))$.
   Similarly, we have
   $$
   eleaf(T)=  |  [k]\setminus \mathcal{I} |+ \sum_{i\in \mathcal{I}}eleaf(T'_i)
   $$
   and
   $$
   plat(\phi(T))=  |  [k]\setminus \mathcal{I} |+ \sum_{i\in \mathcal{I}}plat(\phi(T'_i)).
   $$
   Again by induction hypothesis, we have $eleaf(T)=plat(\phi(T))$, completing the proof.
   \qed

Let $T\in \mathcal{T}_{\mathcal{M}}$. A vertex $i$ of $T$ is said be be {\em even} (resp. {\em odd}) if it is at even (resp. odd) level of $T$.
\begin{theorem}\label{thpsij}
Let $\mathcal{M}=\{1^{k_1}, 2^{k_2}, \ldots, n^{k_n}\}$ with $k_i\geq 1$ for all $1\leq i\leq n$. Suppose that $k_j>1$ for some $j\geq 2$.   Let $M_j$ be the multiset obtained from $\mathcal{M}$ by changing one element  $j$ to   $j-1$, that is,  $\mathcal{M}_j=\{1^{k_1}, 2^{k_2}, \ldots, (j-1)^{k_{j-1}+1}, j^{k_j-1}, \ldots, n^{k_n}\}$. Then there is a bijection $\psi_j$ between   $\mathcal{T}_{\mathcal{M}}$  and $\mathcal{T}_{\mathcal{M}_j}$  such that $$(cdes, casc, eleaf)T= (cdes, casc, eleaf)\psi_j(T)$$ for any $T\in \mathcal{T}_{\mathcal{M}}$.
\end{theorem}

\pf First we give a description of the map $\psi_j$ from  $\mathcal{T}_{\mathcal{M}}$  to $\mathcal{T}_{\mathcal{M}_j}$. Let $T\in \mathcal{T}_{\mathcal{M}}$. In the following,  we demonstrate how to construct $\psi_j(T)$.

\noindent{\bf Case 1.} The odd  vertex $j-1$ is a descendant of the rightmost even vertex  $j$. \\
Suppose that rightmost even vertex $j$ has exactly  $k$ children,  say $x_1, x_2, \ldots, x_k$, listed from left to right.
Let $T_{p}$, $1\leq p\leq k_{j-1}-1$,    be the subtree  rooted at  the $p$-th even vertex $j-1$ (counting from left to right).     Similarly, let $T'_q$, $1\leq q\leq k_{j}-2$,  be the subtree  rooted at the $q$-th even vertex  $j$ (counting from left to right) . First we remove the subtree  $T_{p}$ for all $1\leq p\leq k_{j-1}-1$ and the subtree $T'_q$ for all $1\leq q\leq k_{j}-2$ from $T$.   Then, relabel the rightmost even vertex  $j$ by $j-1$,  and relabel the odd vertex $j$ (resp. $j-1$) by $j-1$ (resp. $j$).  Then we attach  $T_p$ to the odd vertex $j-1$ as its   $p$-th subtree (counting from left to right) for all $1\leq p\leq k_{j-1}-1$, and  attach $T'_q$ to the odd vertex $j$  as its   $q$-th subtree (counting from left to right) for all $1\leq q\leq k_{j}-2$. Let $T^{*}$ be the resulting tree.
In the following, we demonstrate how to get a tree $\psi_j(T)$ from $T^{*}$.

\begin{figure}[h]
\centerline{\includegraphics[width=11cm]{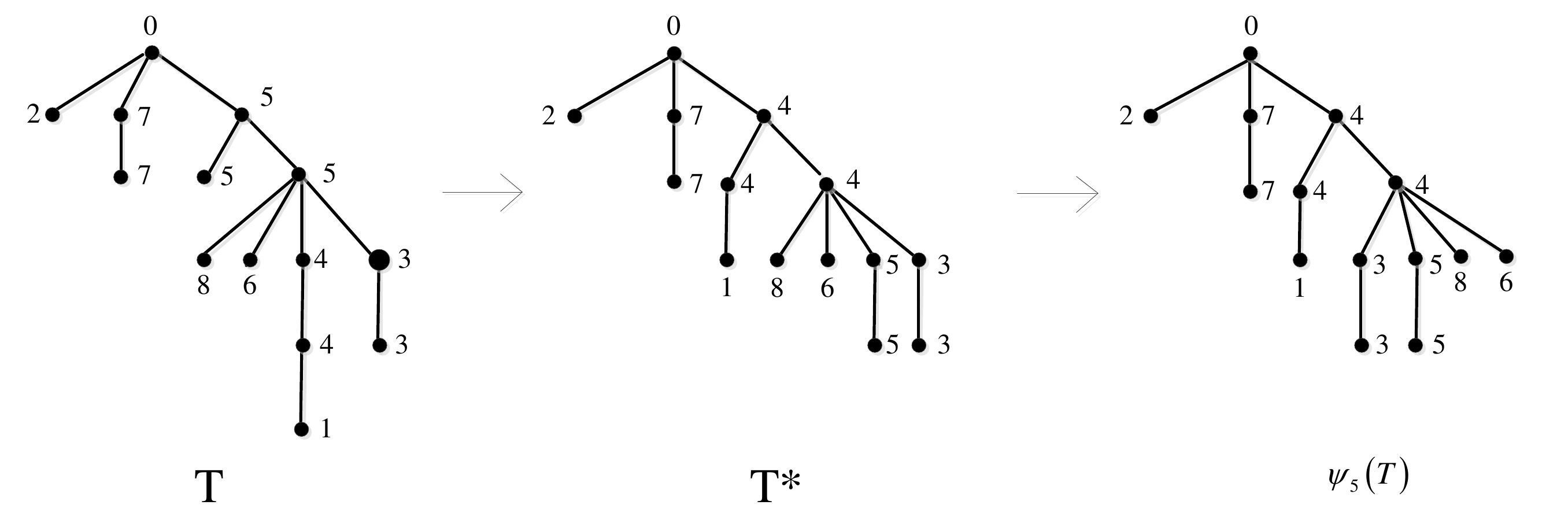}}
\caption{An example of Subcase 1.1.}\label{f1.1}
\end{figure}

 \noindent{\bf Subcase 1.1.} If $x_\ell=j-1$ for some $1\leq \ell\leq k$, then from the construction of $T^{*}$, the odd vertex $j$ is the  $\ell$-th child of the rightmost even vertex $j-1$ in $T^{*}$. Let $\psi_j(T)$ be  the tree obtained from $T^{*}$  by  rearranging the subtrees of the   rightmost even vertex $j-1$ such that the odd vertices $x_{\ell+1}, \ldots,  x_k,  j,   x_1, x_2, \ldots,  x_{\ell-1}$ are its children listed from left to right.   See Figure \ref{f1.1}  for an example.

\begin{figure}[h]
\centerline{\includegraphics[width=11cm]{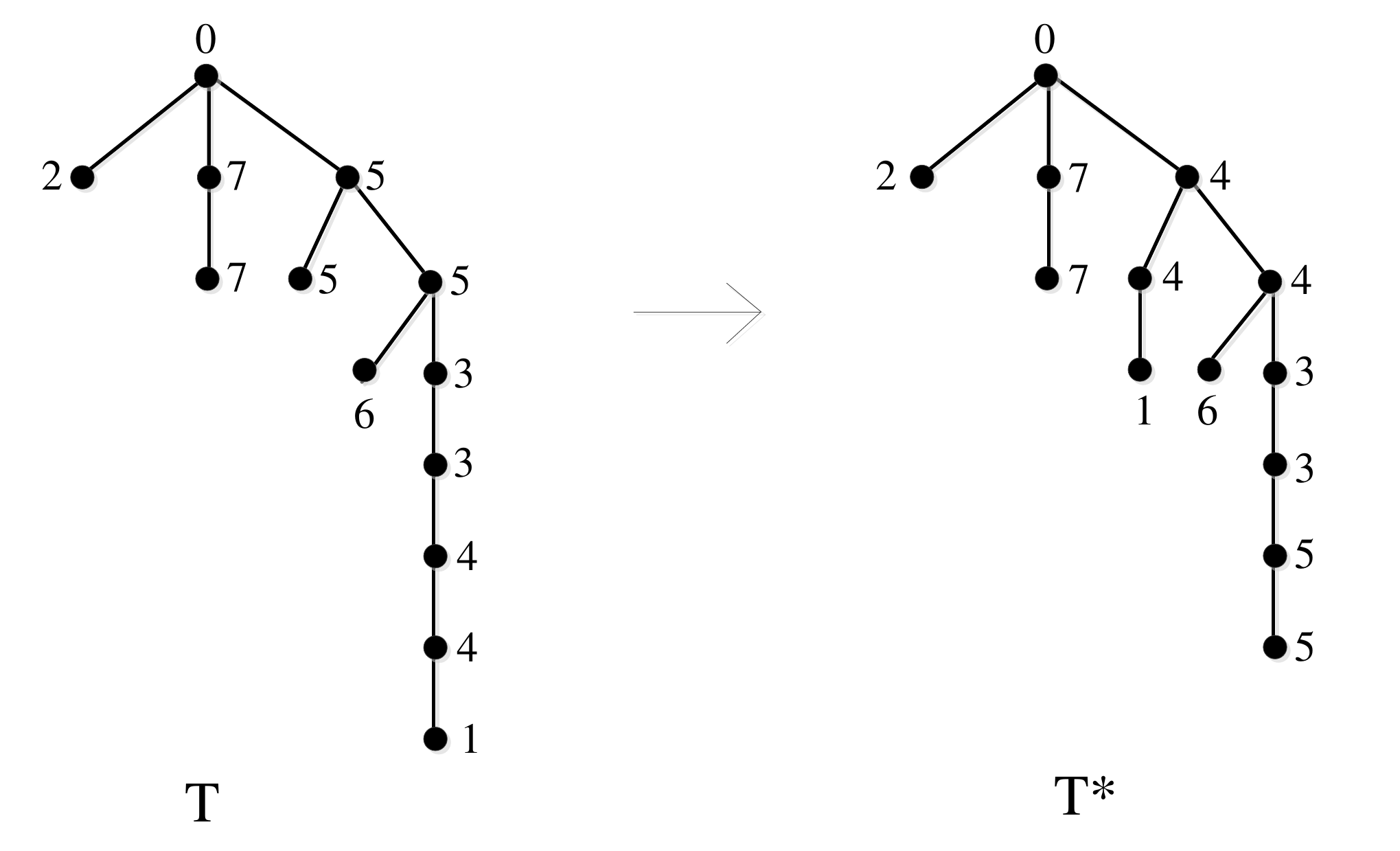}}
 \caption{An example of Subcase 1.2.}\label{f1.2}
\end{figure}

\noindent{\bf Subcase 1.2.}   If $x_{\ell}\neq j-1$ for all $1\leq \ell\leq k$,  set $\psi_j(T)=T^{*}$.  See Figure \ref{f1.2}  for an example.

\begin{figure}[h]
\centerline{\includegraphics[width=11cm]{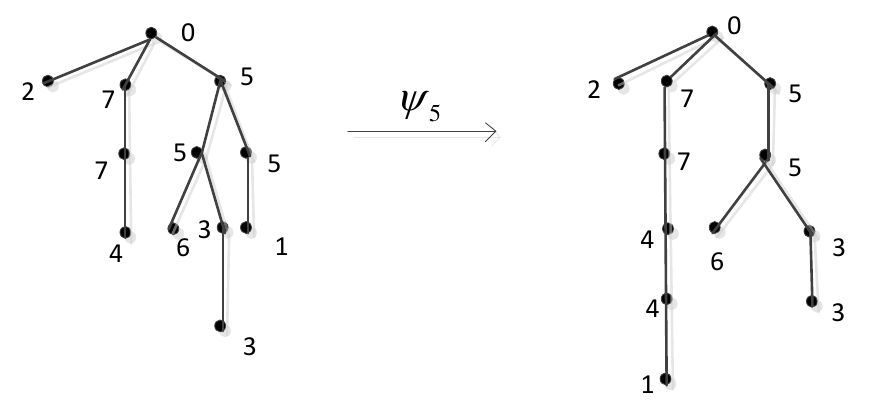}}
\caption{ An example of Case 2. }\label{case2}
\end{figure}

\noindent{\bf Case 2.} Otherwise, let $T'$ be the subtree rooted at rightmost even vertex $j$.
Remove the subtree $T'$ from $T$, relabel the root of $T'$ by $j-1$ and attach it to the odd vertex $j-1$ as its rightmost subtree. Let $\psi_j(T)$ be the resulting tree. See Figure \ref{case2} for an example.

It is routine to check that $\psi_j(T)\in \mathcal{T}_{\mathcal{M}_j}$ and we have $$(cdes, casc, eleaf)T= (cdes, casc, eleaf)\psi_j(T)$$ as desired.

Conversely, given a tree $T\in \mathcal{T}_{\mathcal{M}_j} $, we can recover a tree $T'\in \mathcal{T}_{\mathcal{M}}$ by the following procedure. If the odd  vertex $j$ is a descendant of rightmost even vertex $j-1$, one can  generate a tree $T'\in \mathcal{T}_{\mathcal{M}}$ by reversing the procedure in Case 1. Otherwise, we can recover a tree $T'\in \mathcal{T}_{\mathcal{M}}$ by reversing the procedure in Case 2.  So the construction of the map $\psi_j$ is reversible and  hence the map $\psi_j$ is a bijection. This completes the proof. \qed

\begin{theorem}\label{Psi}
Let $\mathcal{M}=\{1^{k_1}, 2^{k_2}, \ldots, n^{k_n}\}$, $\mathcal{M'}=\{1^{K-n+1}, 2,3,\ldots, n\}$ and $K=k_1+k_2+\ldots+k_n$ with $k_i\geq 1$.
There is a bijection $\Psi$ between $ \mathcal{T}_{\mathcal{M}}$ and $\mathcal{T}_{\mathcal{M}'}$  such that $$(cdes, casc, eleaf)T= (cdes, casc, eleaf)\Psi(T)$$ for any $T\in \mathcal{T}_{\mathcal{M}}$.
    \end{theorem}
\pf  Let $T\in \mathcal{T}_{\mathcal{M}}$. If $k_i=1$ for all $2\leq i\leq n$, then we set $\Psi(T)=T$. Otherwise, find the largest integer $j$ with $j\geq 2$ such that $k_j>2$. By applying the map $\psi_j$ to $T$, we generate a tree $T^{(1)}=\psi_j(T)\in  \mathcal{T}_{\mathcal{M}_j}$, where $M_j$ is the multiset obtained from $\mathcal{M}$ by changing one element  $j$ to   $j-1$. For   $i=1, 2, \ldots, n$, denote by $k'_i$  the number of occurrences of $i$ in $\mathcal{M}_j$.  If $k'_i=1$ for all $2\leq i\leq n$, we stop and set $\Psi(T)=\psi_j(T)$. Otherwise, find the largest integer $\ell$ with $\ell\geq 2$ such that $k'_\ell>2$. By applying the map $\psi_\ell$ to $T^{(1)}$, we generate a tree $T^{(2)}=\psi_\ell(T^{(1)})$. We continue this process until we reach a tree $T'\in T_{\mathcal{M'}}$ and set $\Psi(T)=T'$. By Theorem \ref{thpsij}, the map $\Psi$ is a bijection between between $ \mathcal{T}_{\mathcal{M}}$ and $\mathcal{T}_{\mathcal{M}'}$  satisfying that  $$(cdes, casc, eleaf)T= (cdes, casc, eleaf)\Psi(T)$$ for any $T\in \mathcal{T}_{\mathcal{M}}$ as desired. This completes the proof. \qed

\noindent{\bf Proof of Theorem \ref{mainth1}.} By Theorems \ref{thphi} and \ref{Psi},  the map $\Phi=\phi\circ \Psi\circ \phi^{-1}$ serves as a bijection between $\overline{\mathcal{Q}}_{\mathcal{M}}$ and $\overline{\mathcal{Q}}_{\mathcal{M'}}$. Moreover,  we have $$(asc, des, plat)\pi= (asc, des, plat)\Phi(\pi) $$ for any $\pi\in \overline{\mathcal{Q}}_{\mathcal{M}}$, completing the proof. \qed

In the following, we aim to prove Corollary \ref{coro1}.   The following observation will play an essential role in the proof of Corollary \ref{coro1}.

\begin{observation}\label{ob1}
 For positive integers $m$ and $n$,  let $\mathcal{M}=\{1^{m}, 2,3,\ldots, n\}$. For any permutation $\pi\in \overline{\mathcal{Q}}_{\mathcal{M}}$ with $des(\pi)=n$,  $\pi$ can be uniquely decomposed as
 $$
 \pi=\pi^{(1)}1\pi^{(2)}1\ldots \pi^{(m)}1
 $$
 where each $\pi^{(i)}$ is a (possibly empty)  decreasing sequence.
  \end{observation}

 \noindent{\bf Proof of Corollary \ref{coro1}.}  By Theorem \ref{mainth1}, it suffices to get the enumeration of   permutations $\pi\in  \overline{\mathcal{Q}}_{\mathcal{M}'}$ with $des(\pi)=n$, where $M'=\{1^{K-n+1}, 2,3,\ldots, n\}$. By Observation \ref{ob1}, one can easily verify that the number of
 permutations $\pi\in  \overline{\mathcal{Q}}_{\mathcal{M}'}$ with $des(\pi)=n$ is given by $(K-n+1)^{n-1}$, completing the proof. \qed
\section{Proof of Theorem \ref{mainth2}}
In this section, we shall establish a bijection between quasi-Stirling permutations $\pi\in \overline{\mathcal{Q}}_{\mathcal{M}}$ with $des(\pi)=d+1$ and permutations $\sigma\in  \mathcal{J}_{K,K-n+1}$ with $ exc(\sigma)=d$. To this end, we introduce the {\em  disjoint path and  cycle notation} of a permutation $\sigma\in \mathcal{J}_{n+m-1,m}$. Let $\sigma=\sigma_1\sigma_2\ldots\sigma_{n-1}\in \mathcal{J}_{n+m-1,m} $.  A sequence $(a_1,a_2,\ldots, a_k)$ is called a {\em cycle } of $\sigma$ if $\sigma_{a_i}=a_{i+1}$ for all $1\leq  i\leq k$ with the convention $a_{k+1}=a_1$. A sequence $<a_1,a_2,\ldots a_k>$ is called a {\em path } of $\sigma$ if $\sigma_{a_i}=a_{i+1}$ for all $1\leq  i\leq k-1$  and $a_{k}\geq n$.  Then $\sigma$ can be decomposed into a disjoint union of  distinct paths $P_1, P_2, \ldots, P_s$ and distinct  cycles $C_1, C_2, \ldots, C_t$, written $\sigma=P_1P_2\ldots P_sC_1C_2\ldots C_t$.  Such a decomposition is called the {\em  disjoint path and  cycle notation} of a permutation $\sigma$.  For example, let $\sigma=7816293(11)\in \mathcal{J}_{11, 3}$, the disjoint path and cycle notation of $\sigma$ is given by $<4,6,9><10><5,2,8,11>(1,7,3)$.  Of course this representation of  $\sigma$ in disjoint path and cycle
notation is not unique.  We also have for instance  $\sigma=<10><4,6,9><5,2,8,11>( 7,3,1)$. We can
define a {\em standard representation}   by requiring that (a) the paths are written in increasing order of their  largest element, (b) each cycle  is written with its smallest element first, and (c) the cycles are written in decreasing order of their smallest element.
One can easily check that if $\sigma=P_1P_2\ldots P_sC_1C_2\ldots C_t$ written in standard disjoint path and cycle notation, then we have
\begin{equation}\label{eqexe}
exc(\sigma)=\sum_{i=1}^{s} (asc(P_i)-1) + \sum_{j=1}^{t}  (asc(C_j)-1).
\end{equation}
For instance, if $\sigma=<4,6,9><10><5,2,8,11>(1,7,3)$, then $exc(\sigma)=2+0+2+1=5$.

\begin{lemma}\label{chi}
For positive integers $n$ and $m$,
let $\mathcal{M}=\{1,2,\ldots, n-1, n^{m}\}$. There exists a bijection $\chi$ between quasi-Stirling permutations $\pi\in \overline{\mathcal{Q}}_{\mathcal{M}}$  and permutations $\sigma\in J_{n+m-1, m}$
 such that for any $\pi\in \overline{\mathcal{Q}}_{\mathcal{M}}$, we have $asc(\pi)=exc(\sigma)+1$.
\end{lemma}
\pf First we describe a map $\chi$ from  $ \overline{\mathcal{Q}}_{\mathcal{M}}$  to $ J_{n+m-1, m}$. Let $\pi=\pi_1\pi_2\ldots \pi_{n+m-1}\in \overline{\mathcal{Q}}_{\mathcal{M}} $.  Then $\pi$ can be uniquely decomposed as $\pi=\pi'\pi''$ where there is no $n$ in $\pi''$. Now we can generate a permutation $\chi(\pi)$ in standard disjoint  path and cycle notation by the following procedure.
\begin{itemize}
\item Insert a ``$>$" after each $n$  and insert a ``$<$"  where
appropriate in $\pi'$;
\item Insert a  ``$($" before each left-to-right minimum of $\pi''$ and insert a ``$)$"  where
appropriate;
\item Replace $m$'s occurrences of $n$ by $n, n+1, \ldots, n+m-1$ from left to right.
\end{itemize}
It is easily seen  that the resulting permutation $\chi(\pi)\in \mathcal{J}_{n+m-1, m}$. By (\ref{eqexe}),  we have $asc(\pi)=exc(\chi(\pi))+1$.
For example, let $\pi=46995289173$, we have $\chi(\pi)=<4,6,9><10><5,2,8,11>(1,7,3)$ and $asc(\pi)=exc(\chi(\pi))+1=6$.

In order to show that the map $\chi$ is a bijection, we describe a map $\chi'$ from $ J_{n+m-1, m}$ to $\overline{\mathcal{Q}}_{\mathcal{M}}$.  Let $\sigma\in  J_{n+m-1, m}$ written in standard path and cycle notation. We can generate a permutation $\chi'(\sigma)$ by replacing each  element    larger than $n$  by $n$ and erasing the parentheses and angle brackets. It is easily seen that $\chi'(\sigma)\in \overline{\mathcal{Q}}_{\mathcal{M}} $. Moreover, the standard path and cycle notation  ensures that the maps $\chi$ and $\chi'$ are inverses of each other, and hence the map $\chi$ is a bijection as desired. This completes the proof. \qed

\begin{lemma}\label{delta}
For positive integers $n$ and $m$,
let $\mathcal{M}=\{1^{m},2,\ldots, n\}$. There exists a bijection $\delta$ between quasi-Stirling permutations $\pi\in \overline{\mathcal{Q}}_{\mathcal{M}}$  and permutations $\sigma\in J_{n+m-1, m}$
 such that for any $\pi\in \overline{\mathcal{Q}}_{\mathcal{M}}$, we have $des(\pi)=exc(\sigma)+1$.
\end{lemma}

\pf Let $\pi\in \overline{\mathcal{Q}}_{\mathcal{M}}$ and $M'=\{1,2,\ldots, n^{m}\}$.  First, we generate a permutation $\pi'\in \overline{\mathcal{Q}}_{\mathcal{M}'}$ from $\pi$ by replacing  each element $i$ with $n+1-i$.  It is apparent that we have $des(\pi)=asc(\pi')$. Then,  by applying  the map $\chi$ to $\pi'$, we obtain a permutation $\chi(\pi')\in J_{n+m-1, m} $.  Set $\delta(\pi)=\chi(\pi')$.  By Lemma \ref{chi}, one can easily verify that
  $$des(\pi)=asc(\pi')=exc(\chi(\pi'))+1=exc(\delta(\pi))+1.$$  Since $\chi$ is a bijection, the map $\delta$ is a bijection as desired, completing the proof. \qed

\noindent{\bf Proof of Theorem \ref{mainth2}.}  By Theorem \ref{mainth1} and Lemma \ref{delta}, the map $\delta\circ\Phi$ serves as a bijection between quasi-Stirling permutations $\pi\in \overline{\mathcal{Q}}_{\mathcal{M}}$ with $des(\pi)=d+1$ and permutations $\sigma\in  \mathcal{J}_{K,K-n+1}$ with $ exc(\sigma)=d$. This completes the proof. \qed

\section{Proof of Corollary \ref{coro2}}

For positive integers $m$ and $n$,
let $\mathcal{P}_{m,n}$ denote the set of $m$-tuples $(\pi^{(1)}, \pi^{(2)}, \ldots, \pi^{(m)})$ where
each $\pi^{(i)}$ is a (possibly empty) permutation, $\pi^{(i)}\cap \pi^{(j)}=\emptyset$ for all $1\leq i<j\leq m$ and $\cup_{i=1}^{m}\pi^{(i)}=[n]$.
 For $\alpha=(\pi^{(1)}, \pi^{(2)}, \ldots, \pi^{(m)})\in \mathcal{P}_{m,n}$, let $\mathcal{E}(\alpha)=\{i\mid  \pi^{(i)} =\epsilon\}$.
 Define
 $$
 P_{m,n}(t,u,v)=\sum_{\alpha=(\pi^{(1)}, \pi^{(2)}, \ldots, \pi^{(m)})\in \mathcal{P}_{m,n} }v^{|\mathcal{E}(\alpha)|} \prod_{i\in [m]\setminus \mathcal{E}(\alpha)} t^{des(\pi^{(i)})}u^{asc(\pi^{(i)})}.
 $$
From the definition of $\widetilde{A}(t,u;z)$,  it follows that
 \begin{equation}\label{P}
P_{m,n}(t,u,v)=n![z^n](\widetilde{A}(t,u;z)-1+v)^m.
\end{equation}

 For $1\leq j\leq m$, denote by $\mathcal{P}^{j}_{n,m}$   the set of $m$-tuples $(\pi^{(1)}, \pi^{(2)}, \ldots, \pi^{(m)})\in \mathcal{P}_{n,m}$ such that $\pi^{(j)}$ contains the element $1$. Define
 $$
 P^{j}_{m,n}(t,u,v)=\sum_{\alpha=(\pi^{(1)}, \pi^{(2)}, \ldots, \pi^{(m)})\in \mathcal{P}^j_{m,n} }v^{|\mathcal{E}(\alpha)|} \prod_{i\in [m]\setminus \mathcal{E}(\alpha)} t^{des(\pi^{(i)})}u^{asc(\pi^{(i)})}.
 $$
 It is apparent that
 \begin{equation}\label{eq4.1}
 P^{i}_{m,n}(t,u,v)=P^{j}_{m,n}(t,u,v)
 \end{equation}
 for all $1\leq i<j\leq m$.
 Combining (\ref{P}) and (\ref{eq4.1}), we have
 \begin{equation}\label{eq4.2}
 P^{1}_{m,n}(t,u,v)= {n!\over m}[z^n](\widetilde{A}(t,u;z)-1+v)^m.
 \end{equation}

 \begin{lemma}\label{gamma}
 For positive integers $m$ and $n$, let $\mathcal{M}=\{1^m, 2,\ldots, n\}$.
 There exists a bijection $\zeta$ between $\mathcal{P}^{1}_{m,n}$ and $\overline{\mathcal{Q}}_{\mathcal{M}}$ such that for any $\alpha=(\pi^{(1)}, \pi^{(2)}, \ldots, \pi^{(m)})\in \mathcal{P}^1_{m,n}$, we have $plat(\zeta(\alpha))=|\mathcal{E}(\alpha)|$ and
 $asc(\zeta(\alpha))=\sum_{i\in [m]\setminus \mathcal{E}(\alpha) } asc(\pi^{(i)})$, and $des(\zeta(\alpha))=\sum_{i\in [m]\setminus \mathcal{E}(\alpha) } des(\pi^{(i)})$.
 \end{lemma}
 \pf Let $\alpha=(\pi^{(1)}, \pi^{(2)}, \ldots, \pi^{(m)})\in \mathcal{P}^1_{m,n}$.
  Suppose that  $\pi^{(1)}=\sigma 1\tau$.
  Set $$\zeta(\alpha)=\sigma 1 \pi^{(2)}1 \pi^{(3)}1\ldots \pi^{(m)}1 \tau.$$ It is easily seen that
  $\zeta(\alpha)\in \overline{\mathcal{Q}}_{\mathcal{M}}$ verifying the properties
  \begin{itemize}
  \item $plat(\zeta(\alpha))=|\mathcal{E}(\alpha)|$;
   \item $asc(\zeta(\alpha))=\sum_{i\in [m]\setminus \mathcal{E}(\alpha) } asc(\pi^{(i)})$;
    \item   $des(\zeta(\alpha))=\sum_{i\in [m]\setminus \mathcal{E}(\alpha) } des(\pi^{(i)})$.
   \end{itemize}
   Clearly, the map $\zeta$ is reversible and hence the map $\zeta$ is a bijection, completing the proof. \qed

   \noindent{\bf Proof of Corollary \ref{coro2}.} Let $\mathcal{M'}=\{1^{K-n+1}, 2, \ldots, n\}$.  By Lemma \ref{gamma}, we have
   $$
   \begin{array}{lll}
   \overline{Q}_{\mathcal{M'}}(t,u,v)&=& \sum_{\pi\in \overline{\mathcal{Q}}_{\mathcal{M}'}}t^{des(\pi)}u^{asc(\pi)}v^{plat(\pi)}\\
   &=& \sum_{\alpha=(\pi^{(1)}, \pi^{(2)}, \ldots, \pi^{(K-n+1)})\in \mathcal{P}^1_{K-n+1,n} }v^{|\mathcal{E}(\alpha)|} \prod_{i\in [K-n+1]\setminus \mathcal{E}(\alpha)} t^{des(\pi^{(i)})}u^{asc(\pi^{(i)})}\\
   &=&P^{1}_{K-n+1,n}(t,u,v)\\
   &=& {n!\over K-n+1}[z^n](\widetilde{A}(t,u;z)-1+v)^{K-n+1}  \,\,\,\,\,\, \mbox{(By (\ref{eq4.2}))}.
   \end{array}
   $$
   By Theorem \ref{mainth1}, we have $$\overline{Q}_{\mathcal{M}}(t,u,v)=\overline{Q}_{\mathcal{M'}}(t,u,v)={n!\over K-n+1}[z^n](\widetilde{A}(t,u;z)-1+v)^{K-n+1}$$ as desired. This completes the proof. \qed

 \noindent{\bf Acknowledgments.}
 This work was supported by  the National Natural Science Foundation of China (12071440).


\end{document}